\documentclass[12pt,a4paper]{amsart}
\usepackage{amsfonts,enumerate}
\usepackage{amssymb}
\usepackage{color}

\numberwithin{equation}{section}

\theoremstyle{plain}
\newtheorem{Th}{Theorem}[section]
\newtheorem{Lemma}[Th]{Lemma}

 \theoremstyle{definition}

\newtheorem{Conj}[Th]{Conjecture}
\newtheorem*{Rem}{Remark}
\newtheorem{?}[Th]{Problem}

\newcommand{\Z}{\mathbb{Z}}
\DeclareMathOperator{\supp}{supp}

\newtheorem*{claim}{Claim}
\newtheorem*{cdimo}{Proof of claim}

\newcommand{\bdef}{\begin{definizione}}		
\newcommand{\bteo}{\begin{Th}}			
\newcommand{\eteo}{\end{Th}}			
\newcommand{\bclaim}{\begin{claim}}			
\newcommand{\eclaim}{\end{claim}}			
\newcommand{\bconj}{\begin{Conj}}			
\newcommand{\econj}{\end{Conj}}			
\newcommand{\boss}{\begin{Rem}}			
\newcommand{\eoss}{\end{Rem}}			
\newcommand{\bdimo}{\begin{proof}}			
\newcommand{\edimo}{\end{proof}}			
\newcommand{\bcdimo}{\begin{cdimo}}			
\newcommand{\ecdimo}{\end{cdimo}}			
\newcommand{\blem}{\begin{Lemma}}			
\newcommand{\elem}{\end{Lemma}}			
\DeclareMathOperator{\mmod}{mod}
\allowdisplaybreaks

\begin{document}
\bibliographystyle{amsplain}

\title{Minimizing the frequency of carries in modular addition
}

\author{Francesco Monopoli}       
\address{Dipartimento di Matematica, Universit\`a degli Studi di Milano\\
     via Saldini 50, Milano\\
     I-20133 Italy}

\email{francesco.monopoli@unimi.it}


\begin{abstract}
When adding integers in base $m$, carries occur. The same happens modulo a generic integer $q$ when the set of digits is a complete set of residues modulo $m$ for some positive integer $m$ dividing $q$.
In this paper we prove that asymptotically every digital set in this setting induces carries with frequency at least $1/4$, thus generalizing results of Alon, Diaconis, Shao and Soundararajan.
\end{abstract}
     \maketitle

\section{Introduction}
Fix a nonnegative integer $m$. Adding integers in base $m$ results in carries in the following way: take a set $A$ of $m$ elements, the digits, which is a complete set of residues modulo $m$. We will call such a set a \emph{digital set}.
Adding integers $x_1, x_2$ with last digits $a_1, a_2 \in A$, we find the unique element $a \in A$ such that
    \[ a_1+a_2 \equiv a \pmod m , \] which will be the last digit of the
     sum, and $(a_1 + a_2 -a)/m$ will be the \emph{carry}.

It is natural to look for digital sets which minimize either the number of distinct carries or their frequency at which they generate carries, thus looking for an answer to the following two questions:

{\bf Q1:} What is the minimal number of distinct carries that a digital set $A$ induces, i.e., can we bound from below the quantity
$$C_1(A):=\left|\left\{ \frac{a_1 + a_2 -a}{m} : a_1, a_2, a \in A, a_1 + a_2 \equiv a \mod m\right\}\right|,$$
and what is the structure of digital sets inducing the minimal number of distinct carries?

{\bf Q2:} What is the minimal frequency of carries, i.e., can we bound from below the quantity
$$C_2(A):=\frac{|\{ (a_1, a_2) \in A\times A : a_1 + a_2 \not\in A\}|}{|A|^2},$$
and what is the structure of digital sets inducing the minimal frequency of carries?

It is not hard to prove that for any digital set $A \subseteq \Z$, $C_1(A) \geq 2$ and $m^2C_2(A) \geq \lfloor m^2/4 \rfloor$.
Moreover, up to certain linear transformations, the only digital set which induces only two distinct carries is $[0, m-1]$, whereas the only digital set inducing minimal frequency of carries is $(-m/2, m/2]$ (see \cite{diaconis2014carries}).

The same two questions can be formulated when addition is done modulo a positive integer $q$.
In this setting, for $m|q$, a digital set is a set $A\subseteq \Z_q$, $|A|=m$, which is a complete set of coset representatives for $m\Z_q$.

Diaconis, Shao and Soundararajan \cite{diaconis2014carries}, using an adaptation of the rectification arguments in \cite{BiluLevRuzsa} and \cite{greenruzsa} to this setting, prove that, if $q=p^2$ and $m=p$ for an odd prime $p$, then $C_1(A) \geq 2$ and if $C_1(A)=2$, then there exist $c \in \Z_{p^2}^\times$ and $d\in p\Z_{p^2}$ such that either $cA+d =\{0, 1, \dots, p-1\}$ or $cA+d=\{1, 2, \dots, p\}$, thus giving a complete answer to the first question modulo the square of a prime.

The aforementioned result has been recently generalized to the generic modulus case \cite{monoruzsa}.
When $q$ is not a prime power we have to impose the additional condition on $q$ and $m$, that both are composed of the same primes, and the exponent of each prime in $q$ is strictly greater than in $m$.
This assumption, which does not exclude anyway the case of digital sets of cardinality $m$ modulo $m^2$, is needed in order to avoid trivial cases of digital sets which are either contained in a nontrivial coset, or are union of cosets of a nontrivial subgroup.
Under this hypothesis, the main result in \cite{monoruzsa} states that $C_1(A) \geq 2$, and if $C_1(A)=2$, then there exist $c \in \Z_{q}^\times$ and $d\in m\Z_{q}$ such that either $cA+d =\{0, 1, \dots, m-1\}$ or $cA+d=\{1, 2, \dots, m\}$.

As far as the second question is concerned, Alon \cite{alon} proves that for any digital set $A\subseteq \Z_{p^2}$ of cardinality $p$ for an odd prime $p$, we have $C_2(A) \geq (p^2-1)/4p^2$.

The proof relies on Pollard's inequality for sets $A \subseteq \Z_q$ with the Chowla property, i.e. such that $(a_1 - a_2, q)=1$ for any $a_1, a_2 \in A$.
Digital sets in $\Z_{p^2}$ of cardinality $p$ have this property, but this is no longer true if the cardinality of the digital set is not a prime.

Generalizations of Pollard's inequality for generic modulus exist (see \cite{greenruzsa2}, \cite{grynpol} and \cite{hamiserrapol}), but none of the known bounds is strong enough to be used is an argument similar to the one in \cite{alon} to provide a bound for $C_2(A)$ in the generic modulus case.

The only known bound for $C_2(A)$ in this case can be found in \cite{diaconis2014carries}, where the authors prove that $C_2(A) \geq 2/9$, which is the best bound possible not dependent on the cardinality of the digital set, as can be seen by taking cosets representatives for $3\Z_9 \subseteq \Z_9$.

The main goal of this paper is to prove that the asymptotic bound $C_2(A) \geq 1/4$ for $|A| \rightarrow +\infty$ holds, and this is the best bound possible, as can be seen by taking a symmetric interval as digital set.

\bteo\label{minfreqmodqgen}
 Let $q$ and $m$ be positive integers  composed of the same primes such that the exponent of each prime in $q$ 
  is strictly greater than in $m$.
Let $A \subseteq \Z_q$ be a digital set with $|A|=m$. Let $p^\alpha = \max\{ p_i^{\alpha_i}:  \mbox{$p_i$ prime}, p_i^{\alpha_i} | m\}$ and $\delta_m=1$ if $m$ is odd and $\delta_m=1$ if $m$ is even. Then 
$$C_2(A) \geq \mu(m), $$
where
$$\mu(m)=\begin{cases}
\frac{1 -1/p^{2\alpha}  -2/p^\alpha+\delta_m 2/m}{4} & \mbox{if $p$ is odd,} \\
\frac{1}{4} & \mbox{if $p=2$.}
\end{cases}$$
In particular,
$$\lim_{m \rightarrow + \infty} \min_{|A|=m} C_2(A) =\frac{1}{4}.$$
\eteo

\section{Frequency of carries}
For sets $A, B \subseteq \Z_q$ and $t>0$, let
\begin{equation*}
S(A, B, t)= \sum_{i=1}^t |A+_i B| = \sum_{x \in \Z_q} \min(t, r_{A+B}(x)),
\end{equation*}
where $A+_i B = \{ x \in A+B : r_{A+B}(x) \geq i\}$ is the set of elements in $A+B$ having at least $i$ representations as the sum of an element in $A$ and one in $B$.

Pollard's inequality states that if $A$ or $B$ has the Chowla property, for $1 \leq t \leq \min(|A|, |B|)$ we have
\begin{equation}\label{eqdipol2}
S(A, B, t) \geq t \min(q, |A| + |B| -t).
\end{equation}

As already observed in the introduction, digital sets having nonprome cardinality do not have the Chowla property, and hence \eqref{eqdipol2}, which is the tool used in \cite{alon} to bound $C_2(A)$ in the prime square modulus case, does not necessarly hold.
Nevertheless, in this section we first prove that, if $q=p^\beta$ is a power of a prime, digital sets $A\subseteq \Z_q$, $|A|=p^\alpha$, $0 < \alpha < \beta$, satisfy Pollard's inequality for $t=\lfloor p^\alpha/2 \rfloor$:
\begin{equation}\label{eqdipol3}
S(A, A, t) \geq t(2p^\alpha -t).
\end{equation}
This allows us to prove that for such sets we have
$$C_2(A) \geq \begin{cases} \frac{p^{2\alpha}-1}{4p^{2\alpha}} & \mbox{if $p$ is an odd prime,} \\
\frac{1}{4} & \mbox{if $p=2$.}
\end{cases}$$

In order to do this we need the inverse theorem for Pollard's inequality, due to Nazarewicz, O'Brien, O'Neill and Staples \cite{nazapol}, characterizing pairs of sets $(A, B)$, one of which has the Chowla property, which satisfy equality in \eqref{eqdipol2} for some $t$.

\bteo[Nazarewicz, O'Brien, O'Neill and Staples]\label{nazpolmodq}
Let $A, B \subseteq \Z_q$, $2 \leq t \leq |B| \leq |A|$, $B$ has the Chowla property and $S(A, B, t) = t\min(q, |A|+|B|-t)$. Then one of the following holds:
\begin{enumerate}
\item $|B|=t$,
\item $|A|+|B| \geq q + t $,
\item $|A|=|B|=t+1$, and $B=g-A$ for some $g\in\Z_q$,
\item $A$ and $B$ are arithmetic progression of the same common difference.
\end{enumerate}
\eteo

We remark that the authors in \cite{nazapol} don't prove the theorem in his generality for generic modulus and sets with the Chowla property, but rather restrict their analysis to the case of modulus $q$ prime, so that the Chowla condition in trivially satisfied by any subset of $\Z_q$.
However, it is not a hard task to modify their arguments to prove Theorem \ref{nazpolmodq}.

For odd primes we prove the following:
\bteo\label{4condteo}
Let $A, B \subseteq \Z_{p^\beta}$, $p$ odd, be digital sets, $|A|=|B|=p^\alpha$, with $0 < \alpha < \beta$. Then the following hold:
\begin{enumerate}[(i)]
\item\label{digfatto1} $S\left(A, B, \frac{p^\alpha -1}{2} \right) \geq\frac{3p^{2\alpha} - 2p^\alpha -1}{4}$, 
\item\label{digfatto2}$S\left(A, B, \frac{p^\alpha +1}{2} \right) \geq\frac{3p^{2\alpha} + 2p^\alpha -1}{4}$, 
\item\label{digfatto3} $S\left(A, B, \frac{p^\alpha -1}{2} \right) = \frac{3p^{2\alpha} - 2p^\alpha -1}{4}$ if and only if $A$ and $B$ are arithmetic progressions of the same common difference, 
\item\label{digfatto4}$S\left(A, B, \frac{p^\alpha +1}{2} \right) = \frac{3p^{2\alpha} +2p^\alpha -1}{4}$ if and only if $A$ and $B$ are arithmetic progressions of the same common difference.
\end{enumerate}
\eteo
Before proving the theorem, observe that \eqref{digfatto2}  implies \eqref{digfatto1}, and \eqref{digfatto4} implies \eqref{digfatto3}.
In fact, if \eqref{digfatto2} holds, then
\begin{align*}
S\left(A, B, \frac{p^\alpha -1}{2} \right) &= \sum_{x \in \Z_{p^\beta}} \min\left( \frac{p^\alpha -1}{2}, r_{A+B}(x)\right) \\
&=\sum_{x \in \Z_{p^\beta}} \min\left( \frac{p^\alpha + 1}{2}, r_{A+B}(x)\right) - |A+_{\frac{p^\alpha+1}{2}}B| \\
&\geq \frac{3p^{2\alpha} + 2p^\alpha -1}{4} - p^\alpha \\
&=\frac{3p^{2\alpha} - 2p^\alpha -1}{4},
\end{align*}
where $|A+_{\frac{p^\alpha+1}{2}}B| \leq p^\alpha$ since for every coset $x+\langle p^\alpha \rangle \subseteq \Z_{p^\beta}$ we have
$$\sum_{y \in x+\langle p^\alpha \rangle} r_{A+B}(y) = p^\alpha,$$
and so no more than one element in each coset can have more than $(p^\alpha + 1)/2$ representations in $A+B$.
The same argument shows that  \eqref{digfatto4} implies \eqref{digfatto3}.

To prove Theorem \ref{4condteo} we need some easily verified properties of the $\min$ function, contained in the following lemma.
\begin{Lemma}\label{lemmamin}
Let $a_i,b_i \geq 0, i=1, \dots, n$ and $c \geq 0$. Then
\begin{align}
\min\left(\sum_{i=1}^n a_i, \sum_{i=1}^n b_i\right) \geq \sum_{i=1}^n \min(a_i, b_i) \\
\sum_{i=1}^n \min (c, a_i) \geq \min \left(c, \sum_{i=1}^n a_i\right)
\end{align}
\end{Lemma}
\bdimo
The first inequality is obvious, since the LHS is either $\sum_{i=1}^n  a_i$ or $\sum_{i=1}^n b_i$, and both are clearly greater than $\sum_{i=1}^n \min (a_i, b_i)$.

For the second inequality, if $c \leq a_i$ for some $i$, then
$$\sum_{i=1}^n \min (c, a_i) \geq c = \min \left(c, \sum_{i=1}^n a_i\right),$$
whereas, if $c \geq a_i$ for all $i$, then
$$\sum_{i=1}^n \min (c, a_i) =\sum_{i=1}^n a_i  \geq \min \left(c, \sum_{i=1}^n a_i\right).$$
\edimo

\bdimo[Proof of Theorem \ref{4condteo}]
The proof of Theorem \ref{4condteo} goes by induction on $\alpha$, for all $\beta > \alpha$. For $\alpha=1$ the claims hold by Pollard's inequality and Theorem \ref{nazpolmodq}.

Suppose $\alpha \geq 2$ and let $A_i = i + \langle p \rangle, B_i = i + \langle p \rangle$ for $i \in \Z_p$,
and $A'_i = \frac{A_i - i}{p} \subseteq \Z_{p^{\beta-1}}, B'_i = \frac{B_i - i}{p} \subseteq \Z_{p^{\beta-1}}$.
Then for all $i \in \Z_p$, $|A'_i|=p^{\alpha-1}$, and $A'_i$ is a digital set in $\Z_{p^{\beta-1}}$, for, if two elements $a'_1, a'_2 \in A'_i$ satisfy $a'_1 \equiv a'_2$ modulo $p^{\alpha-1}$, then $a_1 = i + pa'_1 \equiv i + pa'_2 = a_2$ modulo $p^\alpha$, $a_1, a_2 \in A$, a contradiction with the hypothesis of $A$ being a digital set.
The same also holds the $B'_j$'s, and so we can apply the induction hypotheses for $A'_i$ and $B'_j$, $i, j \in \Z_p$.

Let
$$\delta = P(\mbox{$A'_i, B'_j$ are arithmetic progressions of the same common difference})$$ and, for $i, j \in \Z_p$, let
$$\delta_{i,j} = \begin{cases} 1 & \mbox{if $A'_i,B'_j$ are arithmetic progressions of the same common difference,} \\
0 & \mbox{otherwise,}
\end{cases}$$
 so that $\delta p^2 = \sum_{i, j \in \Z_p} \delta_{i,j}$.

Define the map  $$\varphi: \Z_p \times \Z_p \rightarrow \left\{ \frac{p^{\alpha-1}-1}{2}, \frac{p^{\alpha-1}+1}{2} \right\}$$
such that, given $k \in\Z_p$, for exactly $(p-1)/2$ couples $(i, j) \in \Z_p \times \Z_p$ with $i+j \equiv k$ mod $p$ we have $\varphi(i, j) = (p^{\alpha-1}-1)/2$ and for the remaining $(p+1)/2$ couples $(i, j) \in \Z_p \times \Z_p$ with $i+j \equiv k$ mod $p$ we have $\varphi(i, j) = (p^{\alpha-1}+1)/2$. Then, using Lemma \ref{lemmamin} and the inductive hypotheses, we have
\begin{align*}
\sum_{x \in \Z_{p^\beta}} \min\left(\frac{p^\alpha+1}{2}, r_{A+B} (x) \right)  &\geq \sum_{k \in \Z_p} \sum_{x \in k + \langle p \rangle} \sum_{\begin{smallmatrix}
i+j \equiv k \\  \mmod p
 \end{smallmatrix}}  \min\left(\varphi(i,j), r_{A_i+B_j} (x) \right) \\
\geq &\sum_{k \in \Z_p} \sum_{\begin{smallmatrix}
i+j \equiv k \\  \mmod p
 \end{smallmatrix}} \sum_{x \in k + \langle p \rangle} \min\left(\varphi(i,j), r_{A_i+B_j} (x) \right) \\
 \geq& \sum_{k \in \Z_p}\bigg[  \frac{p-1}{2} \frac{p^{\alpha-1}-1}{2}\left( 2p^{\alpha-1} - \frac{p^{\alpha-1}-1}{2}\right) \\
&+ \frac{p+1}{2} \frac{p^{\alpha-1}+1}{2}\left( 2p^{\alpha-1} - \frac{p^{\alpha-1}+1}{2}\right)  + \sum_{\begin{smallmatrix}
i+j \equiv k \\  \mmod p
 \end{smallmatrix}} (1- \delta_{i,j})\bigg] \\
= &\frac{3p^{2\alpha} +2p^\alpha - p^2}{4} + (1-\delta)p^2, 
\end{align*} 
so that if $\delta \leq \frac{3p^2 +1}{4p^2}$ we are done.

For $1 \leq d \leq \frac{p^{\beta-1}-1}{2}$, let $P_{A,d} = \{A'_i : A'_i  \mbox{ is an arithmetic progression of difference $d$}\}$.
Clearly $\sum_d |P_{A,d}| \leq p$, and the same holds for $P_{B,d}$.
Then by Cauchy-Schwarz we have
\begin{align*}\delta p^2 &= \sum_{i \in\Z_p} \sum_{j\in\Z_p} \sum_{d}  \chi((A'_i, B'_j) \in P_{A,d} \times P_{B, d})\\
&= \sum_{d} \sum_{i\in \Z_p} \sum_{j\in \Z_p}   \chi((A'_i, B'_j) \in P_{A,d} \times P_{B, d}) \\
&=\sum_{d} |P_{A,d}||P_{B,d}| \\
&\leq \max_d (|P_{A,d}||P_{B,d}|)^{\frac{1}{2}} \left( \sum_d |P_{A, d}| \sum_d |P_{B, d}| \right)^{\frac{1}{2}} \\
& \leq \max_d (|P_{A,d}||P_{B,d}|)^{\frac{1}{2}} p.\end{align*}
Without loss of generality, after a dilation if necessary we can certainly assume that $\max_d (|P_{A,d}||P_{B,d}|)^{\frac{1}{2}} = (|P_{A,1}||P_{B,1}|)^{\frac{1}{2}}$, since if $A'_i$ is an arithmetic progression of difference $d$ not coprime with $p$, then we would find two elements in $A'_i$ congruent modulo $p^{\alpha-1}$, but this cannot happen in a digital set.

Let $\epsilon_A = |P_{A,1}|/p$ and $\epsilon_B = |P_{B,1}|/p$, $\tilde{P}_A = \{A_i : A'_i \in P_{A, 1} \}, \tilde{P}_B = \{B_j : B'_j \in P_{B, 1} \}$.
Hence
$$\delta \leq \sqrt{\epsilon_A \epsilon_B}
.$$

So, if $\sqrt{\epsilon_A \epsilon_B}\leq \frac{3p^2+1}{4p^2}$ we are done. 

Suppose this does not hold, so that $ \frac{3p^2+1}{4p^2} < \sqrt{\epsilon_A \epsilon_B} \leq\frac{\epsilon_A + \epsilon_B}{2} $.
For $A'_i \in P_{A, 1}, B'_j \in P_{B,1}$ let 
$$A_i = u_i + p \cdot\left[ -\frac{p^{\alpha-1}-1}{2}, \frac{p^{\alpha-1}-1}{2}\right], \qquad B_j = v_j + p \cdot\left[ -\frac{p^{\alpha-1}-1}{2}, \frac{p^{\alpha-1}-1}{2}\right]$$
for some $u_i \in i + \langle p \rangle, v_j \in j + \langle p \rangle$.
Let  
$$U=\{u_i : A'_i \in P_{A,1}\}, \qquad V=\{v_j : B'_j \in P_{B,1}\},$$
and $I_A, I_B$ be the images respectively of $U$ and $V$ under the canonical projection $\pi: \Z_{p^{\beta}} \rightarrow \Z_p$.
For $k\in \Z_p$, let $r_k=r_{I_A + I_B} (k)$.
Since $|I_A|+|I_B| = p(\epsilon_A + \epsilon_B) \geq  p + \frac{p+1}{2}$,
we have $ \frac{p+1}{2}\leq r_k \leq p$ for all $k$. 

Let
$$f(k)=
\frac{p+1}{2}\frac{p^{\alpha-1}+1}{2} + \left(r_k - \frac{p+1}{2}\right)\frac{p^{\alpha-1}-1}{2}=r_k \frac{p^{\alpha-1}-1}{2} + \frac{p+1}{2}. 
$$
Then we can split
\begin{align}\label{dasupereq}
\sum_{x \in \Z_{p^\beta}} \min\left(\frac{p^\alpha+1}{2}, r_{A+B} (x) \right)  \geq& \sum_{k \in \Z_p} \sum_{x \in k + \langle p \rangle} \bigg[ \min\left(f(k), r_{P_{\tilde{A}}+P_{\tilde{B}}} (x) \right) \\
&+  \min \left(
\frac{p^\alpha+1}{2}-f(k), \tilde{r}(x) \right) \bigg], \nonumber
\end{align}
where $\tilde{r}(x)=r_{(A\setminus P_{\tilde{A}}) + P_{\tilde{B}}}(x) + r_{P_{\tilde{A}}+(B\setminus P_{\tilde{B}})}(x) + r_{(A\setminus P_{\tilde{A}}) + (B\setminus P_{\tilde{B}})}(x)$.

Lemmas \ref{lema} and \ref{lema2} below give bounds for the first part of the summation in \eqref{dasupereq}, while Lemma \ref{lema3} provides a bound for the second part.

Using the fact that the representation function of the sum of two intervals $A'_i, B'_j$ is triangular-shaped, we can indeed prove the following:
\blem\label{lema}
Let $k \in \Z_p$ and, with the notation above, let $x_k$ be the element in $k + \langle p \rangle \subseteq \Z_{p^\beta}$ which maximizes $r_{U+V}(x_k)$. 
Then
\begin{align*}\label{triangdiseq} \sum_{x \in k + \langle p \rangle} 
 \min(f(k),  r_{P_{\tilde{A}}+P_{\tilde{B}}} (x))  \geq &r_k \left( \frac{3p^{2\alpha-2}-2p^{\alpha-1}-1}{4}\right)  + \frac{p^\alpha+p^{\alpha-1}}{2} \\
& + \min\left(r_k - \frac{p+1}{2}, r_k - r_{U+V}(x_k)  \right).
\end{align*}
\elem
\bdimo 
Fix $k$
and let $R=r_{U+V}(x_k)$.
Recall that for $A'_i \in P_{A,1}$ and $B'_j \in P_{B,1}$ the representation functions of $A'_i + B'_j$ is triangular-shaped, i.e., it is a translation of the function $\psi(x)=\max(0, p^{\alpha-1}- |x|)$.  
For $y_i \in \Z_{p^{\beta-1}}$ let $\psi_{y_i}(x)=\psi(x-y_i)$.
Then
\begin{equation}\label{quantrk} \sum_{x \in k + \langle p \rangle} \min(f(k), r_{P_{\tilde{A}}+P_{\tilde{B}}} (x)) = \sum_{x \in \Z_{p^{\beta-1}}} \min\left(f(k), \sum_{\begin{smallmatrix}
i\in I_A, j \in I_B \\
i+j \equiv k \\  \mmod p
\end{smallmatrix}} r_{A'_i+B'_j} (x)\right).\end{equation}

To get the desired bound for \eqref{quantrk} we minimize 
\begin{equation}\label{quantrk2}S(\mathcal{P}_k) =  \sum_{x \in \Z_{p^{\beta-1}}} \min\left(f(k),\sum_{i=1}^{r_k} \psi_{y_i}(x) \right),\end{equation}
where $\mathcal{P}_k$ ranges over all possible multisets of $r_k$ elements $y_1, \dots, y_{r_k} \in \Z_{p^{\beta-1}}$ with the condition that $R$ $y_i$'s are equal.
Without loss of generality, up to a translation, we can assume that those $R$ points are equal to $0$.

First of all, we compute $S(\bar{\mathcal{P}}_{k,a,b})$ for
$$\bar{\mathcal{P}}_{k,a,b} =\{\underbrace{-1, \dots, -1}_{a}, \underbrace{0, \dots, 0}_{R }, \underbrace{1, \dots, 1}_{b}\},$$
with the conditions that $a+R+b=r_k$ and $a+R, b+R \geq r_k/2$, which imply that $a, b \leq r_k/2$.
In this case, we have
$$\sum_{i=1}^{r_k} \psi_{y_i}(x)= \begin{cases}
a & \mbox{$x=-p^{\alpha-1}$}, \\
a(p^{\alpha-1}+x+1)+ R(p^{\alpha-1}+x) + b(p^{\alpha-1}+x-1) & \mbox{$x\in [-p^{\alpha-1}+1, -1]$}, \\
a(p^{\alpha-1}-1) + Rp^{\alpha-1} + b(p^{\alpha-1}-1) & \mbox{$x=0$}, \\
a(p^{\alpha-1}-x-1) + R(p^{\alpha-1}-x) + b(p^{\alpha-1}-x+1) & \mbox{$x \in [1, p^{\alpha-1}-1]$}, \\
b & \mbox{$x=p^{\alpha-1}$}.\\
\end{cases}$$
For $x \in [1, p^{\alpha-1}-1]$ we have $\sum_{i=1}^{r_k} \psi_{y_i}(x) \geq f(k)$ if and only if
$$1\leq x \leq \frac{p^{\alpha-1}+1}{2} + \frac{2(b-a)-p-1}{2r_k}.$$
We have $ \frac{2(b-a)-p-1}{2r_k} \in (-2, 0)$, since $b-a \leq b \leq r_k/2 \leq \frac{p-1}{2} < \frac{p}{2}$ and
$r_k \geq \frac{p+1}{2}$.
\\
{\bf Case 1}: $\frac{2(b-a)-p-1}{2r_k} \in (-2, -1)$.

Since $-\frac{p+1}{2r_k} \in [-1, 0)$, this can happen only if $a > b$.
In this case, for $x \in [1, p^{\alpha-1}-1]$, we have $\sum_{i=1}^{r_k} \psi_{y_i}(x) \geq f(k)$ if and only if $x \leq \frac{p^{\alpha-1}-3}{2}$.
Hence
\begin{align*}
\sum_{x \in [1, p^{\alpha-1}]}\min\left(f(k),  \sum_{i=1}^{r_k} \psi_{y_i}(x) \right)  =& f(k) \frac{p^{\alpha-1}-3}{2} + \sum_{x=\frac{p^{\alpha-1}-1}{2}}^{p^{\alpha-1}-1} (-r_k x + p^{\alpha-1}r_k + b -a) + b \\
=&f(k)\frac{p^{\alpha-1}-3}{2} + \frac{p^{\alpha-1}+1}{2}(b-a+p^{\alpha-1}r_k)\\
& - \frac{r_k}{2}\frac{3p^{2\alpha-2}-3}{4} + b.
\end{align*}
\\
{\bf Case 2}: $\frac{2(b-a)-p-1}{2r_k} \in [-1, 0)$.

In this case, for $x \in [1, p^{\alpha-1}-1]$, we have $\sum_{i=1}^{r_k} \psi_{y_i}(x) \geq f(k)$ if and only if $x \leq \frac{p^{\alpha-1}-1}{2}$.
Hence
\begin{align*}
\sum_{x \in [1, p^{\alpha-1}]}\min\left(f(k),  \sum_{i=1}^{r_k} \psi_{y_i}(x) \right) = &f(k) \frac{p^{\alpha-1}-1}{2} + \sum_{x=\frac{p^{\alpha-1}+1}{2}}^{p^{\alpha-1}-1} (-r_k x + p^{\alpha-1}r_k + b -a) + b \\
=&f(k)\frac{p^{\alpha-1}-1}{2} + \frac{p^{\alpha-1}-1}{2}(b-a+p^{\alpha-1}r_k) \\
&- \frac{r_k}{2}\frac{3p^{2\alpha-2}-4p^{\alpha-1}+1}{4} + b.
\end{align*}

The same argument gives a similar computation for $\sum_{x \in [-p^{\alpha-1}, -1]} \sum_{i=1}^{r_k} \psi_{y_i}(x)$, where the roles of $a$ and $b$ are exchanged.

Since we cannot fall into Case $1$ both for $[-p^{\alpha-1},-1]$  and $[1,p^{\alpha-1}]$, a simple computation shows that, if we fall into Case $1$ for $[1,p^{\alpha-1}]$ and Case $2$ for $[-p^{\alpha-1},-1]$, we have 
\begin{align*}
\sum_{x \in [-p^{\alpha-1}, p^{\alpha-1}]}\min\left(f(k),  \sum_{i=1}^{r_k} \psi_{y_i}(x) \right) =&r_k \left( \frac{3p^{2\alpha-2}-2p^{\alpha-1}-1}{4}\right) + \frac{p^\alpha+p^{\alpha-1}}{2} \\
&+ r_k - \frac{p+1}{2} + 2b \\
\geq& r_k \left( \frac{3p^{2\alpha-2}-2p^{\alpha-1}-1}{4}\right) + \frac{p^\alpha+p^{\alpha-1}}{2} \\
&+ r_k - \frac{p+1}{2},
\end{align*}
whereas, if we fall in Case $2$ both for $[1,p^{\alpha-1}]$ and $[-p^{\alpha-1},-1]$, then
\begin{align*}
\sum_{x \in [-p^{\alpha-1}, p^{\alpha-1}]}\min\left(f(k),  \sum_{i=1}^{r_k} \psi_{y_i}(x) \right) =&r_k \left( \frac{3p^{2\alpha-2}-2p^{\alpha-1}-1}{4}\right) + \frac{p^{\alpha}+p^{\alpha-1}}{2} \\
&+ a+b \\
=& r_k \left( \frac{3p^{2\alpha-2}-2p^{\alpha-1}-1}{4}\right) + \frac{p^{\alpha}+p^{\alpha-1}}{2} \\
& + r_k - R,
\end{align*}
thus proving the lemma in this particular case.

We now show that for any multiset $\mathcal{P}_k=\{y_1, \dots, y_{r_k}\}$ with $R$ equal elements, we have
$S(\mathcal{P}_k) \geq S(\bar{\mathcal{P}}_{k, a, b})$
for some choices of $a,b$ with $a+b+R=r_k$.

Observe first of all that if 
$\bar{x} \in \Z_{p^{\beta-1}}$
 satisfies $\sum_{i=1}^{r_k} \psi_{y_i}(\bar{x})  >  f(k)$, then

\begin{align*}
f(k)=r_k\frac{p^{\alpha-1}-1}{2} + \frac{p+1}{2} &< \sum_{\begin{smallmatrix} 
y \in \mathcal{P}_k: \\
\bar{x} \in \supp \psi_{y}
\end{smallmatrix}} \psi_y (\bar{x})\\
& \leq p^{\alpha-1} \cdot |\{ y \in \mathcal{P}_k : \bar{x} \in \supp(\psi_y)\}|,
\end{align*}
so that
$$|\{ y \in \mathcal{P}_k : \bar{x} \in \supp(\psi_y)\}| \geq \frac{r_k}{2} + \frac{p+1-r_k}{2p^{\alpha-1}}> \left\lfloor \frac{r_k}{2} \right\rfloor.$$
Since $\{ y \in \mathcal{P}_k : \bar{x} \in \supp(\psi_y)\} \subseteq [\bar{x} -( p^{\alpha-1}-1), \bar{x}+ (p^{\alpha-1}-1)]$, this implies that if $x\in \Z_{p^{\beta-1}}$ satisfies $\sum_{i=1}^{r_k} \psi_{y_i}(x)  >  f(k)$, then $x \in [\bar{x} -2( p^{\alpha-1}-1), \bar{x}+ 2(p^{\alpha-1}-1)]$.
Hence, if $\mathcal{P}_k=\{y_1, \dots, y_{r_k} \}$ is such that $S(\mathcal{P}_k)$ is minimal, then we can assume that $\supp \left( \sum_{i=1}^{r_k} \psi_{y_i} \right) \subseteq [\bar{x} - 4(p^{\alpha-1}-1), \bar{x} + 4(p^{\alpha-1}-1)]$, for if this does not happen, we can take any $y \not\in  [\bar{x} -3( p^{\alpha-1}-1), \bar{x}+ 3(p^{\alpha-1}-1)]$ and replace it with an element in $[\bar{x} -3( p^{\alpha-1}-1), \bar{x}+ 3(p^{\alpha-1}-1)]$ to obtain a multiset $\mathcal{P}'_k$ with $S(\mathcal{P}'_k) \leq S(\mathcal{P}_k)$.

Note that $ [\bar{x} - 4(p^{\alpha-1}-1), \bar{x} + 4(p^{\alpha-1}-1)] \subsetneq \Z_{p^{\beta-1}}$ with the exception of a finite number of choices of $(p, \alpha, \beta)$, where the conclusions of the lemma can be easily checked, and hence we have that
either there exists no element $x \in \Z_{p^{\beta-1}}$ with $\sum_{i=1}^{r_k} \psi_{y_i}(x) > f(k)$, and so we get the trivial bound $r_k p^{2\alpha-2}$ in \eqref{quantrk} and the lemma is true, or we can assume $\supp \left( \sum_{i=1}^{r_k} \psi_{y_i} \right) \subsetneq \Z_{p^{\beta-1}}$.

Suppose we are in the latter case, so that we can assume $\supp \left( \sum_{i=1}^{r_k} \psi_{y_i} \right) \subseteq \mathcal{I}$, where $\mathcal{I}$ is an interval different from the whole $\Z_{p^{\beta-1}}$, and so for $x_1, x_2 \in  \mathcal{I}$, $x_1 < x_2$ has the obvious meaning.

Let $\sigma(x)=\sum_{i=1}^{r_k} \psi_{y_i}(x)$.
Consider the minimal element $y_1 \in \mathcal{P}_k$, say with multeplicity $m_1$, and consider the multiset $\mathcal{P}'_k $ obtained from $\mathcal{P}_k$ by translating $y_1$ to $y_1+1$, i.e.,
$\mathcal{P}'_k = \{ y'_i\}$,
with 
$$y'_i = \begin{cases} y_i +1& \mbox{if $y_i = y_1$,}\\
 y_i & \mbox{otherwise.}
\end{cases}$$
We will show that  $S(\mathcal{P}'_k) \leq S(\mathcal{P}_k)$, which is equivalent, for $\sigma'(x)=\sum_{i=1}^{r_k} \psi_{y'_i}(x)$, to
$$\sum_{x \in \Z_{p^{\beta-1}}} \max(\sigma'(x)-f(k), 0) -  \max(\sigma(x)-f(k), 0) \geq 0.$$
If $\sigma(x) < f(k)$ for all $x \leq y_1$ the claim holds, since then for all $x\in \Z_{p^{\beta-1}}$  $\sigma(x) < f(k)$ or $\sigma'(x) \geq \sigma(x)$.

Let $\bar{t}$ be the largest nonnegative integer such that $\sigma(y_1-\bar{t}) > f(k)$.
If $\bar{t}=0$, the claim holds since $\sigma(y_1+1) +m_1 \geq \sigma(y_1)$.

Let $\bar{t} > 0$ and $M= |\{y \in \mathcal{P}_k: y \neq y_1, y_1 - \bar{t}\in \supp(\psi_y) \}|$.
Then the following hold:
\begin{enumerate}
\item $\sigma(y_1+t) \geq \sigma(y_1-t) + 2M$ for $t \in [1, \bar{t}+1]$,
\item $\sigma(y_1 - t) \geq \sigma(y_1-t-1) + m_1 + M$ for $t\in [0, \bar{t}]$,
\item $\sigma(y_1-\bar{t}) = \sigma(y_1 - \bar{t}-1)+M+m_1$.
\end{enumerate}
{\bf Case 1:} $\sigma'(y_1-\bar{t})=\sigma(y_1-\bar{t}) - m_1 \geq f(k)$.

In this case we have $\sigma(y_1 + \bar{t}+1) \geq \sigma(y_1 - \bar{t}-1)+2M = \sigma(y_1 - \bar{t}) -m_1+ M \geq f(k)$.
Hence for all $x \in[y_1 - \bar{t}, y_1 + \bar{t} +1]$ we have $\sigma(x), \sigma'(x) \geq f(k)$ and so
$$
\sum_{x \in[y_1 - \bar{t}, y_1 + \bar{t} +1]} \max(0,  \sigma'(x)-f(k)) - \max(0, \sigma(x)-f(k))= \sum_{x \in [y_1 - \bar{t}, y_1]}\!\!\!\!\! -m_1 + \sum_{x \in (y_1, y_1 +\bar{t} + 1]}\!\!\!\!\! m_1 =0.
$$
\\
{\bf Case 2:} $\sigma'(y_1-\bar{t})=\sigma(y_1-\bar{t}) - m_1 < f(k)$.

Since, arguing as above, $\sigma(y_1 + \bar{t}+1) +m_1 \geq \sigma(y_1 - \bar{t}) \geq f(k)$, we have
\begin{align*}
\sum_{x \in[y_1 - \bar{t}, y_1 + \bar{t} +1]} \max(0, &\, \sigma'(x)-f(k)) - \max(0, \sigma(x)-f(k))=-\sigma(y_1 - \bar{t})+f(k) \\
&+ \sum_{x \in [y_1 - \bar{t} +1, y_1]} -m_1 + \sum_{x \in (y_1, y_1 + \bar{t}]} m_1  \\
&+ \sigma(y_1+\bar{t}+1) +m_1 -f(k) - \tau \geq 0,
\end{align*}
where 
$$\tau = \begin{cases} 0 & \mbox{if } \sigma(y_1 + \bar{t}+1) < f(k), \\
			\sigma(y_1+\bar{t} +1) - f(k) & \mbox{otherwise.}
\end{cases}$$

In both cases we have $S(\mathcal{P}'_k) \leq S(\mathcal{P}_k)$ since for all $x \not \in [y_1 - \bar{t}, y_1 + \bar{t} +1]$  $\sigma(x) < f(k)$ or $\sigma'(x) \geq \sigma(x)$.
%

Suppose that $R > r_k/2$.
Say we have $a$ elements $y \in \mathcal{P}_k$, $y < 0$ and $b$ elements $y \in \mathcal{P}_k, y >0$.
Then, iterating the shifting procedure explained above, which has an obvious equivalent for the maximal element of $\mathcal{P}_k$, if we replace those $a$ elements with $y' = -1$ and those $b$ elements with $y' = 1$, we recover a translate of the multiset $\bar{\mathcal{P}}_{k,a,b}$ with $a+b=r_k-R$ and we have
$S(\bar{\mathcal{P}}_{k,a,b}) \leq S(\mathcal{P}_k)$, so that the conclusion of the lemma holds.

Suppose that $R \leq  r_k/2$.
Then all $y \in \mathcal{P}_k$ have multiplicity less or equal to $r_k/2$, and, with the shifting procedure explained above, we iteratively shift the minimal element $y_1$ of $\mathcal{P}_k$, and we stop doing this as soon as we get $y_1+1=y_2$, $y_1<y_2$ with multiplicity respectively $m_1$ and $m_2$ and $m_1+m_2\leq r_k/2$.
We do the same thing for the maximal element of $\mathcal{P}_k$, and we end up with a new multiset $\mathcal{P}'_k$ with $S(\mathcal{P}'_k) \leq S(\mathcal{P}_k)$ having at most three distinct elements, each with multiplicity less or equal to $r_k/2$.
Indeed, since the sum of the first two consecutive distinct elements in $\mathcal{P}'_k$ is strictly greater than $r_k/2$, the sum of the multiplicities of the remaining elements must be less or equal to $r_k/2$, and so, if there is more than one of these elements, we could certainly shift the maximal element at least one more time to the left.

This new multiset $\mathcal{P}'_k$ is equal, up to a translation, to $\bar{\mathcal{P}}_{k,a,b}$ for some $a, b$ satisfying $a, b \leq r_k/2$ and $R'= r_k -a-b \leq r_k/2$.

Since 
$$\min\left(r_k - \frac{p+1}{2}, r_k -R'\right) = r_k-\frac{p+1}{2} = \min\left(r_k - \frac{p+1}{2}, r_k -R\right),$$
the conclusion follows from the computation at the beginning of the proof.
\edimo

To sum the contributions given by Lemma \ref{lema}, we need the following:
\blem\label{lema2}
\begin{equation}\label{eqinlema2}\sum_{k \in \Z_{p}} \min\left(r_k - \frac{p+1}{2}, r_k - r_{U+V}(x_k)  \right) \geq \frac{p+1}{2}\left((\epsilon_A + \epsilon_B)p-\frac{3p+1}{2}\right).\end{equation}
\elem
\bdimo
Let $\{R_1, \dots, R_l\} = \{r_{U+V}(x_{k}) : r_{U+V}(x_{k}) > (p+1)/2 \}$.

Since $r_{U+V}(x) \leq p$ for all $x$, we have that $R_1, \dots, R_l$ are the $l$ highest values among $\{r_{U+V}(x): x \in \Z_{p^\beta}\}$.


Since $U$ and $V$ have the Chowla property, we have
\begin{align*}
|U||V| - \frac{p+1}{2}\left( |U|+|V| - \frac{p+1}{2}\right) & \geq |U||V| - \sum_{i=1}^{\frac{p+1}{2}}  |U+_{i} V| \\
&= \sum_{i=\frac{p+3}{2}}^{\min(|U|, |V|)}  |U+_{i} V| \\
& = \sum_{x \in \Z_{p^\beta}} \max\left(r_{U+V} (x) - \frac{p+1}{2}, 0\right) \\
&=- l\frac{p+1}{2} + \sum_{k=1}^l R_k 
\end{align*}

Recalling that $\sum_{k \in \Z_p}r_k = |U||V| = \epsilon_A \epsilon_B p^2$, we have
\begin{align*}
\sum_{k \in \Z_{p}} \min\left(r_k - \frac{p+1}{2}, r_k - r_{U+V}(x_k)  \right)=& \sum_{k: r_{U+V}(x_k) \leq \frac{p+1}{2}} \left( r_k - \frac{p+1}{2}\right)  \\
&+ \sum_{k: r_{U+V}(x_k) > \frac{p+1}{2}} \left( r_k - r_{U+V}(x_k) \right) \\
\geq & \sum_{k \in \Z_p} r_k -(p-l)\frac{p+1}{2} -l \frac{p+1}{2} - \epsilon_A \epsilon_B p^2 \\
& + \frac{p+1}{2} \left( (\epsilon_A + \epsilon_B)p - \frac{p+1}{2}\right) \\
=&\frac{p+1}{2}\left( (\epsilon_A + \epsilon_B)p-\frac{3p+1}{2}\right)
\end{align*}
as required.
\edimo

\blem\label{lema3}
\begin{align*}\sum_{k \in \Z_p} \sum_{x \in k + \langle p \rangle} \min\left(\frac{p^\alpha+1}{2}-f(k), \tilde{r}(x) \right) & \geq \epsilon_A (1-\epsilon_B)p^2 \bigg( \frac{3p^{2\alpha-2}-2p^{\alpha-1}-1}{4}+1\bigg) \\
&+  (1-\epsilon_A)\epsilon_Bp^2 \bigg( \frac{3p^{2\alpha-2}-2p^{\alpha-1}-1}{4}+1\bigg) \\
&+(1-\epsilon_A)(1-\epsilon_B)p^2 \bigg( \frac{3p^{2\alpha-2}-2p^{\alpha-1}-1}{4}\bigg) 
\end{align*}
\elem
\bdimo
Since for every $k$ we have $\frac{p^\alpha+1}{2}-f(k) = (p-r_k)\left( \frac{p^{\alpha-1} -1}{2}\right)$, we can compute
\begin{align*}\sum_{k \in \Z_p} \sum_{x \in k + \langle p \rangle} \min\left(\frac{p^\alpha+1}{2}-f(k), \tilde{r}(x) \right)  \geq \sum_{k \in \Z_p} & \sum_{\begin{smallmatrix}
i\not\in I_A, j \in I_B \\
i+j \equiv k \\  \mmod p
\end{smallmatrix}} \min\left(\frac{p^{\alpha-1}-1}{2}, r_{A'_i+B'_j}(x) \right) \\
&+\sum_{\begin{smallmatrix}
i\in I_A, j \not\in I_B \\
i+j \equiv k \\  \mmod p
\end{smallmatrix}} \min\left(\frac{p^{\alpha-1}-1}{2}, r_{A'_i+B'_j}(x) \right) \\
& + \sum_{\begin{smallmatrix}
i\not\in I_A, j \not\in I_B \\
i+j \equiv k \\  \mmod p
\end{smallmatrix}} \min\left(\frac{p^{\alpha-1}-1}{2}, r_{A'_i+B'_j}(x) \right).
\end{align*}
Using the inductive hypothesis to get a bound better then the one coming from Pollard's inequality whenever we consider the sumset $A'_i + B'_j$ with $i \not\in I_A, j \in I_B$ or viceversa, we get the desired bound.
\edimo

Using Lemmas \ref{lema2} and \ref{lema3}
we can finish the proof of the Theorem \ref{4condteo}:

\begin{equation}\label{ultimaeq}
\begin{aligned}
\sum_{x \in \Z_{p^\beta}} \min\bigg(\frac{p^\alpha+1}{2}, &\,  r_{A+B} (x) \bigg)  \geq \sum_{k \in \Z_p} \sum_{x \in k + \langle p \rangle} \bigg[ \min\left(f(k), r_{P_{\tilde{A}}+P_{\tilde{B}}} (x) \right) \\
&+\min \left(\frac{p^\alpha+1}{2}-f(k), \tilde{r}(x)\right)  \bigg] \\
\geq&  \epsilon_A \epsilon_B p^2 \bigg( \frac{3p^{2\alpha-2}-2p^{\alpha-1}-1}{4}\bigg)  + p\frac{p^\alpha+p^{\alpha-1}}{2} \\
& + \frac{p+1}{2}\left( (\epsilon_A + \epsilon_B)p-\frac{3p+1}{2}\right) \\
&+ \epsilon_A (1-\epsilon_B)p^2 \bigg( \frac{3p^{2\alpha-2}-2p^{\alpha-1}-1}{4}+1\bigg) \\
&+  (1-\epsilon_A)\epsilon_Bp^2 \bigg( \frac{3p^{2\alpha-2}-2p^{\alpha-1}-1}{4}+1\bigg) \\
&+(1-\epsilon_A)(1-\epsilon_B)p^2 \bigg( \frac{3p^{2\alpha-2}-2p^{\alpha-1}-1}{4}\bigg) \\
\geq&  \bigg( \frac{3p^{2\alpha}+2p^{\alpha}-1}{4}\bigg) + p^2\bigg(\frac{\epsilon_A +\epsilon_B}{2} + \epsilon_A(1-\epsilon_B) + \epsilon_B(1-\epsilon_A) -1\bigg) \\
&+ p\bigg( \frac{\epsilon_A+\epsilon_B}{2}-1\bigg).
\end{aligned}
\end{equation}
Since
\begin{align*} p^2\bigg(\frac{\epsilon_A +\epsilon_B}{2} + \epsilon_A(1-\epsilon_B) + \epsilon_B(1-\epsilon_A) -1\bigg) 
+ p\bigg( \frac{\epsilon_A+\epsilon_B}{2}-1\bigg) \\
\geq p^2\bigg(-2\epsilon_A \epsilon_B + 3\sqrt{\epsilon_A \epsilon_B}) - 1\bigg) 
+ p\bigg( \sqrt{\epsilon_A \epsilon_B}-1\bigg),
\end{align*}
and $x^2(-2p^2) + x(3p^2+p) -p^2-p \geq 0$
for $1 \geq x \geq (p+1)/2p$, the conclusion holds since we assumed $\sqrt{\epsilon_A\epsilon_B} > (3p^2+1)/4p^2 \geq (p+1)/2p$.

This concludes the proof of \eqref{digfatto2}.

To prove \eqref{digfatto4}, notice that, in order to have equality in \eqref{digfatto4}, from \eqref{ultimaeq} we must have $\epsilon_A\epsilon_B=1$, so that every $A'_i$ and $B'_j$ is an arithmetic progression of the same common difference $d$, $(d, p)=1$, and after a dilation if necessary we can assume $d=1$.
Moreover, by Theorem \ref{nazpolmodq}, since we must also have equality in \eqref{eqinlema2}, we have that both $U$ and $V$ are arithmetic progressions of the same common difference $d'$, $(d', p)=1$, say $U=\{u_0, u_1 =u_0+d', \dots, u_{p-1}=u_0+(p-1)d'\}$ and $V=\{v_0, v_1=v_0+d', \dots, v_{p-1}=v_0+(p-1)d'\}$.

From the proof of part \eqref{digfatto2}, since $u_0+v_0 \equiv u_1 + v_{p-1}$ modulo $p$, we deduce that $pd' \equiv \pm p$ modulo $p^\beta$, so that $d' \equiv \pm 1$ modulo $p^{\beta-1}$, and $A$ is an arithmetic progressions of difference $d'$, starting either from $u_0 - \frac{p^\alpha-p}{2}$ or $u_{0} + \frac{p^\alpha-p}{2}$, and the same holds for $B$, thus completing the proof of the theorem. 
\edimo

We are left to study the case of $p=2$, which is way easier and forms the following theorem:
\bteo\label{teo4condeven}
Let $A, B \subseteq \Z_{2^\beta}$ be digital sets, $|A|=|B|=2^\alpha$, with $0 < \alpha < \beta$. Then
\begin{enumerate}[(i)]
\item $S\left(A, B, 2^{\alpha-1} \right) \geq 2^{2\alpha} - 2^{2\alpha-2},$
\item $S\left(A, B, 2^{\alpha-1} \right) = 2^{2\alpha} - 2^{2\alpha-2}$  if and only if $A$ and $B$ are arithmetic progressions of the same common difference.
\end{enumerate}
\eteo
\bdimo
The proof goes by induction of $\alpha$, for all $\beta>\alpha$.
For $\alpha=1$ the claim holds.

Suppose $\alpha \geq 2$.
\\
Let
$$A_0=A \cap \langle 2 \rangle, \qquad A_1 = A \cap (1+\langle 2 \rangle),$$
$$B_0=B \cap \langle 2 \rangle, \qquad B_1 = B \cap (1+\langle 2 \rangle).$$
Then, $A'_i = \frac{A_i-i}{2} \subseteq \Z_{2^{\beta-1}}, B'_j = \frac{B_j-i}{2} \subseteq \Z_{2^{\beta-1}}$ are digital sets of cardinality $2^{\alpha-1}$ in $\Z_{2^{\beta-1}}$, and thanks to the induction hypothesis we have
\begin{align*}
\sum_{x \in \Z_{2^\beta}} \min\left(2^{\alpha-1}, r_{A+B} (x) \right)  \geq&  \sum_{x \in \langle 2 \rangle} \min\left(2^{\alpha-2}, r_{A_0+B_0} (x) \right) + \min\left(2^{\alpha-2}, r_{A_1+B_1} (x) \right) \\
&+\sum_{x \in 1+ \langle 2 \rangle} \min\left(2^{\alpha-2}, r_{A_0+B_1} (x) \right) + \min\left(2^{\alpha-2}, r_{A_1+B_0} (x) \right)  \\
\geq & 4 \cdot 2^{\alpha-2} (2^\alpha - 2^{\alpha-2}) \\
=& 2^{2\alpha} - 2^{2\alpha-2}
\end{align*} 
as required.

Moreover, by the induction hypothesis, equality holds in the chain of inequalities above if and only if $A'_0, A'_1, B'_0, B'_1$ are arithmetic progressions of the same common difference.
A simple analysis shows that the only possibility for $A$ and $B$ to achieve the equality $S\left(A, B, 2^{\alpha-1} \right) = 2^{2\alpha} - 2^{2\alpha-2}$ under this additional condition is that both $A$ and $B$ are arithmetic progression of the same common difference as claimed.
\edimo

Arguing as in \cite{alon}, Theorems \ref{4condteo} and \ref{teo4condeven} allow us to prove the exact bound for $C_2(A)$ in the case of digital sets $A \subseteq \Z_{p^\beta}$ of cardinality $p^{\alpha}$, proving that such sets induce at least $\left\lfloor p^\alpha/4 \right\rfloor$ carries, with equality holding if and only if $A$ is a dilation of $[-(p^{\alpha-1}-1)/2, (p^{\alpha-1}-1)/2]$ by a factor $d$ coprime with $p$ if $p$ is odd, or a dilation of $[-2^{\alpha-1}, 2^{\alpha -1})$ or $(-2^{\alpha-1}, 2^{\alpha -1}]$ by an odd factor $d$ if $p=2$.

We can now prove Theorem \ref{minfreqmodqgen}.
\bdimo[Proof of Theorem \ref{minfreqmodqgen}.]
Let $m=m' p^\alpha$ for $p^\alpha = \max\{ p_i^{\alpha_i}:  \mbox{$p_i$ prime}, p_i^{\alpha_i} | m\}$, and $A_i = A \cap i + \langle m' \rangle$ for $i=0, \dots, m'-1$.
Writing $A= \{ a_j\}_{j=0, \dots m-1}$, where $a_j \equiv j$ mod $m$ for all $j=0, \dots, m-1$, then for all $i$, $A_i = \{ a_{i}, a_{i+m'}, \dots, a_{i+(p^\alpha-1)m'}\}$, and so $|A_i|=p^\alpha$.
Then $A'_i := \frac{A_i - i}{m'} \subseteq \Z_{q/m'}$, since $A_i -i \subseteq m'\Z_q \simeq \Z_{q/m'}$.
Moreover, for $x, y \in A'_i$,  we have $x \not \equiv y$ modulo $p^\alpha$, for otherwise we would have $i+m'x, i+m'y \in A$ with $i+m'x \equiv i+m'y$ mod $m$, which contradicts the fact that $A$ is a digital set.
Hence $A'_i\subseteq \Z_{q/m'}$ is a digital set for every $i$.

Consider the projection $\pi: \Z_{q/m'} \rightarrow \Z_{p^{\beta}}$, where $p^\beta$ is the highest power of $p$ dividing $q/m'$, $\beta > \alpha$.

We have that $|A'_i| = p^\alpha = |\pi(A'_i)|$, and still for $x, y \in A'_i$, we have $\pi(x) \not\equiv \pi(y)$ modulo $p^\alpha$, so that $\pi(A'_i)$ is, once again, a digital set for every $i$.
\\
{\bf Case 1}: $p$ odd.

Using Theorem \ref{4condteo} and Lemma \ref{lemmamin}, we have
\begin{align*}
\sum_{x \in \Z_q} \min\left(\frac{p^\alpha-1}{2}, r_{A_i+A_j} (x) \right) & = \sum_{x \in \Z_{q/m'}} \min\left(\frac{p^\alpha-1}{2}, r_{A'_i+A'_j} (x) \right) \\
& = \sum_{y \in \Z_{p^\beta}} \sum_{x \in \pi^{-1} (y)} \min\left(\frac{p^\alpha-1}{2}, r_{A'_i+A'_j} (x) \right) \\
& \geq \sum_{y \in \Z_{p^\beta}} \min\left(\frac{p^\alpha-1}{2}, \sum_{x \in \pi^{-1} (y)} r_{A'_i+A'_j} (x) \right) \\
& =  \sum_{y \in \Z_{p^\beta}} \min\left(\frac{p^\alpha-1}{2}, r_{\pi(A'_i)+\pi(A'_j)} (y) \right)  \\
& \geq  \frac{3p^{2\alpha}-2p^\alpha-1}{4}.
\end{align*} 

Using this inequality and Lemma \ref{lemmamin}, we have
\begin{align*}
\sum_{x \in \Z_q} \min\left(\left\lfloor\frac{m}{2}\right\rfloor, r_{A+A} (x) \right) & = \sum_{k=0}^{m'-1} \sum_{x \in k + \langle m' \rangle} \min\left(\left\lfloor\frac{m}{2}\right\rfloor, r_{A+A} (x) \right)  \\
&=  \sum_{k=0}^{m'-1}\sum_{x \in k + \langle m' \rangle} 
 \min \left( \left\lfloor\frac{m}{2}\right\rfloor,\sum_{\begin{smallmatrix}
i+j \equiv k \\  \mmod m'
 \end{smallmatrix}} r_{A_i + A_j}(x)\right)  \\
& \geq  \sum_{k=0}^{m'-1}\sum_{x \in k + \langle m' \rangle} \sum_{\begin{smallmatrix}
i+j \equiv k \\  \mmod m'
 \end{smallmatrix}}
 \min \left( \frac{p^\alpha-1}{2}, r_{A_i + A_j}(x)\right)  \\
& \geq m'^2 \frac{3p^{2\alpha}-2p^\alpha-1}{4} = \frac{3m^2 - 2m^2/p^\alpha - m'^2}{4},
\end{align*}
so that
\begin{align*}\frac{3m^2 - 2m^2/p^\alpha - m^2/p^{2\alpha}}{4} & \leq \sum_{x \in A+A} \min\left(\left\lfloor\frac{m}{2}\right\rfloor, r_{A+A}(x)\right) \\
& \leq \sum_{x \in (A+A) \cap A} \left\lfloor\frac{m}{2}\right\rfloor + \sum_{x \in (A+A) \setminus A} r_{A+A}(x)
\end{align*}
and thus
$$\sum_{x \in (A+A) \setminus A} r_{A+A}(x) \geq m^2\frac{1 -1/p^{2\alpha}  -2/p^\alpha+\delta_m 2/m}{4},$$
where $\delta_m=1$ if $m$ is odd and $\delta_m=0$ if $m$ is even.\\
{\bf Case 2}: $p=2$.

Using Theorem \ref{teo4condeven}, we have
\begin{align*}
\sum_{x \in \Z_q} \min\left(2^{\alpha-1}, r_{A_i+A_j} (x) \right) & = \sum_{x \in \Z_{q/m'}} \min\left(2^{\alpha-1}, r_{A'_i+A'_j} (x) \right) \\
& = \sum_{y \in \Z_{2^{\alpha+1}}} \sum_{x \in \pi^{-1} (y)} \min\left(2^{\alpha-1}, r_{A'_i+A'_j} (x) \right) \\
& \geq \sum_{y \in \Z_{2^\beta}} \min\left(2^{\alpha-1}, \sum_{x \in \pi^{-1} (y)} r_{A'_i+A'_j} (x) \right) \\
& =  \sum_{y \in \Z_{2^\beta}} \min\left(2^{\alpha-1}, r_{\pi(A'_i)+\pi(A'_j)} (y) \right)  \\
& \geq  2^{2\alpha} - 2^{2\alpha-2}.
\end{align*} 

Using this inequality and Lemma \ref{lemmamin}, we get
\begin{align*}
\sum_{x \in \Z_q} \min\left(\frac{m}{2}, r_{A+A} (x) \right) & = \sum_{k=0}^{m'-1} \sum_{x \in k + \langle m' \rangle} \min\left(\frac{m}{2}, r_{A+A} (x) \right)  \\
&=  \sum_{k=0}^{m'-1}\sum_{x \in k + \langle m' \rangle} 
 \min \left( \frac{m}{2},\sum_{\begin{smallmatrix}
i+j \equiv k \\  \mmod m'
 \end{smallmatrix}} r_{A_i + A_j}(x)\right)  \\
& \geq  \sum_{k=0}^{m'-1}\sum_{x \in k + \langle m' \rangle} \sum_{\begin{smallmatrix}
i+j \equiv k \\  \mmod m'
 \end{smallmatrix}}
 \min \left( 2^{\alpha-1}, r_{A_i + A_j}(x)\right)  \\
& m'^2 (2^{2\alpha} - 2^{2\alpha-2})=\frac{3m^2}{4}.
\end{align*}
Then we get
\begin{align*}\frac{3m^2}{4} &\leq \sum_{i=1}^{\frac{m}{2}} |A+_iA|  \\
& = \sum_{x \in A+A} \min\left(\frac{m}{2}, r_{A+A}(x)\right) \\
& \leq \sum_{x \in (A+A) \cap A} \frac{m}{2} + \sum_{x \in (A+A) \setminus A} r_{A+A}(x)
\end{align*}
and so
$$\sum_{x \in (A+A) \setminus A} r_{A+A}(x) \geq \frac{m^2}{4}.$$

Since $\sum_{x \in (A+A) \setminus A} r_{A+A}(x)$ counts the couples $(a_1, a_2) \in A\times A$ such that $a_1+a_2 \not\in A$, i.e. the number of occurrences of carries induced by $A$,
we get the desired conclusion.

For an integer $m$ let $\varphi(m) = \max\{ p_i^{\alpha_i}:  \mbox{$p_i$ prime}, p_i^{\alpha_i} | m\}$ be the largest prime power dividing $m$, $\psi(m)$ be the largest prime dividing $m$ and $\omega(m)$ be the function counting the number of distinct prime factors in $m$.
It's easy to see that $\varphi(m) \rightarrow \infty$ as $m\rightarrow \infty$.
In fact, suppose this does not hold, and let $\{m_i\}$ be an increasing sequence of integers with $\varphi(m_i) \leq L$ for all $i$.
Then, if $\omega(m_i) \rightarrow \infty$, so does $\psi(m_i)$, and consequently the same holds for $\varphi(m_i)$, a contradiction.
On the other hand, if, up to subsequences, $\omega(m_i) \leq M$ for all $m_i$, then clearly $m_i \leq \omega(m_i)\varphi(m_i) \leq LM$, which does not go to infinity, thus leading to a contradiction.

Then, since balanced digital sets of cardinality $m$ induce $\left\lfloor m^2/4\right\rfloor$ carries and the function $\mu(m)$ tends to $1/4$ as $m$ goes to infinity, we have
$$\left\lfloor\frac{m^2}{4} \right\rfloor \frac{1}{m^2} \geq\min_{|A|=m} C_2(A) 
\geq \mu(m)
\longrightarrow \frac{1}{4} \qquad \mbox{for $m$} \longrightarrow \infty,$$
thus completing the proof of the theorem.
\edimo
%
%
%
\section*{}
\bibliography{bibliography}

\providecommand{\bysame}{\leavevmode\hbox to3em{\hrulefill}\thinspace}
\providecommand{\MR}{\relax\ifhmode\unskip\space\fi MR }
\providecommand{\MRhref}[2]{%
  \href{http://www.ams.org/mathscinet-getitem?mr=#1}{#2}
}
\providecommand{\href}[2]{#2}
\begin{thebibliography}{1}

\bibitem{alon}
N.~Alon, \emph{Minimizing the number of carries in addition}, SIAM J. Discrete
  Math. \textbf{27} (2013), no.~1, 562--566. \MR{3035466}

\bibitem{BiluLevRuzsa}
Y.~F. Bilu, V.~F. Lev, and I.~Z. Ruzsa, \emph{Rectification principles in
  additive number theory}, Discrete Comput. Geom. \textbf{19} (1998), no.~3,
  Special Issue, 343--353, Dedicated to the memory of Paul Erd{\H{o}}s.
  \MR{1608875 (2000a:11018)}

\bibitem{diaconis2014carries}
P.~Diaconis, X.~Shao, and K.~Soundararajan, \emph{Carries, group theory, and
  additive combinatorics}, The American Mathematical Monthly \textbf{121}
  (2014), no.~8, 674--688.

\bibitem{greenruzsa2}
B.~Green and I.~Z. Ruzsa, \emph{Sum-free sets in abelian groups}, Israel J.
  Math. \textbf{147} (2005), 157--188. \MR{2166359 (2006e:11030)}

\bibitem{greenruzsa}
\bysame, \emph{Sets with small sumset and rectification}, Bull. London Math.
  Soc. \textbf{38} (2006), no.~1, 43--52. \MR{2201602 (2006i:11027)}

\bibitem{grynpol}
D.~J. Grynkiewicz, \emph{On extending {P}ollard's theorem for
  {$t$}-representable sums}, Israel J. Math. \textbf{177} (2010), 413--439.
  \MR{2684428 (2011k:11016)}

\bibitem{hamiserrapol}
Y.~O. {Hamidoune} and O.~{Serra}, \emph{{A note on Pollard's Theorem}}, arXiv
  (2008).

\bibitem{monoruzsa}
F.~{Monopoli} and I.~Z. {Ruzsa}, \emph{{Carries and the arithmetic progression
  structure of sets}}, arXiv (2015).

\bibitem{nazapol}
E.~Nazarewicz, M.~O'Brien, M.~O'Neill, and C.~Staples, \emph{Equality in
  {P}ollard's theorem on set addition of congruence classes}, Acta Arith.
  \textbf{127} (2007), no.~1, 1--15. \MR{2289969 (2008e:11014)}

\end{thebibliography}
     \end{document}